\documentclass[draft]{publmathd}
\usepackage{amsmath,amsfonts,amssymb}

\pyear{2004}
\pvol{1}

\refno{3285}

 \pagespan{1}{13}
 \shortlefthead{V. Bovdi, C. H\"ofert \es W. Kimmerle}
 \shortrighthead{On the first Zassenhaus conjecture}
 \receivedinfo{July 3, 2004; revised October 22, 2004}

\newtheorem*{theore}{Theorem}
\newtheorem{Theorem}{Theorem}

\newtheorem{Proposition}{Proposition}

\theoremstyle{definition}

\newtheorem*{rema}{Remark}

\def\gp#1{\langle #1 \rangle}

\allowdisplaybreaks
\providecommand{\bdshs}{\belowdisplayshortskip=10pt}

\title[On the first Zassenhaus conjecture for integral group rings]{On the first Zassenhaus conjecture for integral group
rings}

\author[V. Bovdi, C. H\"ofert \es W. Kimmerle]
{V. BOVDI (Debrecen), C. H\"OFERT (Stuttgart) \\ \es W. KIMMERLE (Stuttgart)}

\address{Victor Bovdi\\
    Institute of Mathematics\\
    University of Debrecen\\
    H-4010 Debrecen, P.O. Box 12\\
    Hungary\\
    Institute of Mathematics and Informatics\\
    College of Ny\'\i regyh\'aza\\
    S\'ost\'oi \'ut 31/b, H-4410 Ny\'\i regyh\'aza\\
    Hungary}

\email{vbovdi@math.klte.hu}

\address{C.~H\"ofert\\
    IGT, Fachbereich Mathematik\\
    Universit\"at of Stuttgart\\
    Pfaffenwaldring 57, D-70550, Stuttgart\\
    Germany}
\email{hoefercn@mathematik.uni-stuttgart.de}

\address{W.~Kimmerle\\
    IGT, Fachbereich Mathematik\\
    Universit\"at of Stuttgart\\
    Pfaffenwaldring 57, D-70550, Stuttgart\\
    Germany}
\email{kimmerle@mathematik.uni-stuttgart.de}

\dedicatory{Dedicated to the memory of Professor B\'ela Brindza}

\thanks{The research was supported by OTKA  No.\ T 037202, No.\ T 038059 and partially by the Deutsche Forschungsgemeinschaft}


\begin{document}

\begin{abstract}
It was conjectured by H.~Zassenhaus that a torsion unit of an
integral group ring of a finite group is conjugate to a group
element within the rational group algebra. The object of this note
is the computational aspect of a method developed by I.~S.~Luthar
and I.~B.~S.~Passi which sometimes permits an answer to this
conjecture. We illustrate the method on certain explicit examples.
We prove with additional arguments that the conjecture is valid
for any 3-dimensional crystallographic point group. Finally we
apply the method to generic character tables and establish a  $p$-variation of the conjecture for the simple groups $PSL(2,p)$.
\end{abstract}

\maketitle

\section{Introduction}

Let $V(\mathbb ZG)$ be the group of units of augmentation 1 of an
integral group ring $\mathbb ZG$ of a  finite group $G$. With
respect to the structure of $V(\mathbb ZG)$ H.~Zassenhaus stated
around 1976 the following  conjectures.

\begin{itemize}
\item[]
\begin{itemize}
\item[{\bf ZC--1.}]
Let $u$ be a unit of finite order of $V(\mathbb ZG)$. Then $u$ is conjugate within $\mathbb QG$ to an element of $G$.\footnote{Elements of $G$ are called the trivial units of $\mathbb Z G $.}
\item[{\bf ZC--2.}]
Let $H$ be a subgroup  of $V(\mathbb ZG)$ with the
same order as $G$. Then $H$ is conjugate within $\mathbb QG$ to $G$.
\item[{\bf ZC--3.}]
Let $H$ be a finite subgroup  of $V(\mathbb ZG)$. Then $H$ is conjugate within $\mathbb QG$ to a subgroup of $G$.
\end{itemize}
\end{itemize}
It was shown by K. W.~Roggenkamp and L. L.~Scott that the conjectures ZC--2 and ZC--3 are in general not true \cite{15}, \cite{18}, \cite{12}.

No counterexample however is known to the conjecture ZC--1. The conjectures ZC--2 and ZC--3 are valid\footnote{We say that ZC--i is valid for a
finite group $G$ if $V(\mathbb ZG)$ has the property stated in ZC--i.} for many important classes of finite groups, e.g.\ nilpotent groups \cite{16}, \cite{21}. The results indicate that the following variation of the Zassenhaus conjectures -- a Sylow-like theorem for $V(\mathbb ZG)$ -- may be true.

\smallskip
\noindent
{\bf  p-ZC--3.} Let $H$ be a finite $p$-subgroup  of $V(\mathbb ZG)$. Then $H$ is conjugate \\
   $\quad \phantom{\text{\bf p-ZC--3.}\ }$ within $\mathbb Q G$ to a $p$-subgroup of $G$.
\smallskip

For a recent survey on the Zassenhaus conjectures and variations of them we refer to \cite{11}.

\smallskip
One of our goals is to provide algorithms and programs for the GAP -- package LAGUNA which decide for a given group $G$ whether the conjecture ZC -- i is valid or not \cite{10}. For computational aspects in integral group rings related to the Zassenhaus conjectures see \cite{3}, \cite{2}.

The main purpose of this note is to discuss the method developed by {\sc I. S. Luthar} and {\sc I. B. S. Passi} for the conjecture ZC--1 \cite{13} under the view of computational aspects. We illustrate the Luthar--Passi method  for $A_4 \times S_3 $ and for the octahedral group $ S_4 \times C_2 $. In the second case the method
alone does not suffice to establish the Zassenhaus conjecture
whereas in the first case it suffices.

In \cite{9} it is shown that ZC--1 is valid for all finite groups of order less than 71. We remark that the Luthar--Passi method plays an essential role to establish this result.

It seems to be unknown whether ZC--1 is true for a direct product of finite groups $G \times H$ provided it holds for each of its factors. We show that this is the case when $H = C_2 $.
As a  consequence, we obtain that the conjecture ZC--1 is valid for all finite 3-dimensional
crystallographic point groups.

The method may also be applied to generic character tables. In \cite{20} the simple groups $PSL(2,p^f)$ have been studied with it. In the last section it is shown that for the linear groups $PSL(2,p) $ the variation $p$-ZC--3 holds for the describing characteristic $p $.

\section{The Luthar--Passi method}

The ingredients for this method are the following results.

\begin{Theorem}[\cite{13}]\label{A}
 Suppose that $u$ is an element of\/\ $V(\mathbb ZG)$ of order~$k$. Let $z$ be a primitive $k$-th root  of
unity. Then for every integer $l$ and any ordinary character $\chi$ of $G$,
the number
\begin{equation}\label{e:1}
\mu_l(u,\chi)=\frac{1}{k} \sum_{d|k}Tr_{\mathbb Q(z^d)/\mathbb Q}
\{\chi(u^d)\cdot z^{-dl}\}
\end{equation}
is a non-negative integer. Let $\varphi $ be a
$\mathbb C$ - representation which affords $\chi$. Then
$\mu_l(u,\chi)$ is the multiplicity of $z^l$ in the Jordan
canonical form of $\varphi(u)$. In particular the degree of $\chi$ bounds $\mu_l(u,\chi)$.
\end{Theorem}

Let $u=\sum\alpha_gg $ be a normalized non-central torsion unit of order $k$
and let $\nu_i=\varepsilon_{C_i}(u)$ be the partial augmentation\footnote{Let $u=\sum\alpha_gg \in \mathbb Z G $. Then the partial augmentation
with respect to the conjugacy class $C_i$ is defined as
$\nu_i=\varepsilon_{C_i}(u) = \sum_{g \in C_i} \alpha_g $. }
of $u$ with respect to the conjugacy class $C_i $. Then by well known
theorems of {\sc G. Higman} and {\sc S.~D.~Berman}  \cite[Theorem 10, p.\ 102]{17}
$$
\nu_1 = 0  \quad \text{and more general} \quad \nu_j = 0 $$
if the class $C_j$ consists of a central element. Because $u$ is normalized this implies
      $$ \nu_2+\nu_3+\cdots+\nu_m=1, $$
where $m$ denotes the class number of $G$.

\begin{Theorem}[\cite{5}]\label{B}
Let $u$ be a torsion unit of\/\ $V(\mathbb ZG)$.  The order of $u$ divides the exponent of\/\ $G$.
\end{Theorem}

\begin{Theorem}[{\cite[Theorem 2.7]{14}}]\label{C}
 Let $u$ be a torsion unit of\/\ $V(\mathbb ZG)$.
Let $C$ be a conjugacy class of $G$. If $p$ is a prime
dividing  the order of a representative of $C$ but not the
order of $u$ then the partial augmentation $\varepsilon_C(u)$ is zero.
\end{Theorem}

Now the key result in order to establish the conjecture ZC--1 is the following one.

\begin{Theorem}[\cite{13}, {\cite[Theorem 2.5]{14}}]\label{D}
Let $u$ be a normalized unit of $\mathbb Z G$ of order $k$. Then $u$ is conjugate
in $\mathbb Q G$ to an element $g \in G$ if and only if for each $d$
dividing $k$ there is precisely one conjugacy class $C_{i_d}$ with
partial augmentation $\varepsilon_{C_{i_d}}(u^d)  \neq 0 $.
\end{Theorem}

Starting with the ordinary character table of a finite group $G$ Theorem A yields restrictions on the multiplicities
$\mu_l(u,\chi) $. These multiplicities are via
$$ \bdshs
\chi(u^d) =  \sum_j \varepsilon_{C_j}(u^d) \chi(g_j), \qquad g_j \in C_j $$
related with the partial augmentations of the conjugacy classes. The
restrictions on the multiplicities lead to bounds on the partial
augmentations. Additional information on the partial augmentation
comes inductively from the quotients of $G $  and theoretical
statements like Theorem~B. The starting point of the induction
is given by the fact that ZC--1 is valid for the nilpotent
quotients of $G$ \cite{21}. For some groups, e.g.\ $A_4 \times S_3$
this leads finally via Theorem C to ZC--1 or to statements for
elements of $p$-power order which are relevant with respect to the
conjecture p-ZC--3 (see the last section). The Luthar--Passi method is also of interest in the context of constructing a possible counterexample to ZC--1. It gives precise information about the
partial augmentation of a candidate for a counterexample.

\vspace{-10pt}

\section{The Group $A_4 \times S_3$}

In this section we use the Luthar--Passi method to establish the
following result.

\begin{Proposition}\label{Prop:1}
Every normalized torsion unit  $u$ in $\mathbb ZG$, is rational conjugate to a group element, where $G=A_4\times S_3$.
\end{Proposition}
\def\fvon{\vrule height 2.4ex depth.8ex width0ex}
\def\tvon{\vrule height 2.2ex depth.7ex width0ex}

\begin{proof}
The group $G$ has the following character table (easily verified with \cite{8}):

\medskip
\small{\centerline{ \vbox {{ \offinterlineskip \halign { $#$ & \tvon\vrule \quad
\hfil  $#$ &  \vrule   \quad \hfil  $#$ & \vrule \quad \hfil  $#$ &
\vrule   \quad \hfil  $#$  & \vrule \quad \hfil  $#$ &  \vrule
\quad \hfil  $#$ & \vrule \quad \hfil $#$ &  \vrule   \quad \hfil
$#$ & \vrule \quad \hfil $#$ & \vrule   \quad \hfil $#$& \vrule
\quad \hfil $#$ & \vrule   \quad \hfil $#$\cr \noalign {\hrule}
&1a &3a &3b &2a &2b & 6a & 6b &2c &3c & 3d&  3e &6c\fvon\cr \noalign {\hrule}
\chi_1    &1 & 1 & 1 & 1 & 1 &  1 &  1 & 1 & 1 &  1 &  1 & 1 \cr
\chi_2    &1 & 1 & 1 & 1 &-1 & -1 & -1& -1 & 1 &  1 &  1 & 1 \cr
\chi_3    &1 & \omega& \bar{\omega}  &1  &1 &  \omega & \bar{\omega}
& 1 & 1 &  \omega & \bar{\omega}&  1 \cr %
\chi_4    &1 & \omega& \bar{\omega} & 1 &-1&  -\omega
&-\bar{\omega}& -1 & 1&   \omega & \bar{\omega} & 1 \cr
\chi_5    &1 &\bar{\omega} & \omega & 1 & 1&  \bar{\omega} &  \omega
& 1&  1 & \bar{\omega} &  \omega & 1 \cr %
\chi_6     &1 &\bar{\omega} & \omega & 1& -1& -\bar{\omega} &
-\omega &-1 & 1 & \bar{\omega} & \omega & 1\cr %
\chi_7      & 2 & 2 & 2 & 2 & 0&0&0&0& -1 & -1 & -1 &-1\cr
\chi_8    & 2 & 2\omega &\bar{2\omega} & 2 & 0 &  0 &  0 & 0 &-1&
-\omega& -\bar{\omega} &-1 \cr %
\chi_9    &2 &\bar{2\omega}&  2\omega&  2&  0 &  0 &  0 & 0 &-1
&-\bar{\omega} & -\omega &-1\cr
\chi_{10}    &3 & 0 & 0& -1& -3&   0 &  0 & 1 & 3 &  0 &  0 &-1\cr
 \chi_{11}    &3 & 0 & 0& -1 & 3 &  0&   0& -1 & 3 &
0 &  0& -1\cr  \chi_{12}    & 6 & 0 & 0 &-2 & 0&   0
&  0 & 0 &-3 &  0 &  0&  1\cr
\noalign {\hrule} }}} }}%

\normalsize
\bigskip
\noindent
where $\omega$ denotes a primitive third root of unity.

Note that by Theorem B the possible orders of non-trivial torsion units of
$V(\mathbb Z [A_4 \times S_3])$ are $2$, $3$ or $6$, and  according to Theorem C
we get:
\begin{align}
&\nu_{3a}=\nu_{3b}=\nu_{3c}=\nu_{3d}=\nu_{3e}=\nu_{6x}=0 \quad \text{when}\  k=2; \nonumber\\
&\nu_{2a}=\nu_{2b}=\nu_{2c}=\nu_{6x}=0 \quad  \text{when}\  k=3;\label{e:2}\\
&\nu_{3a}=\nu_{3b}=\nu_{2a}=\nu_{2b}=\nu_{2c}=\nu_{3c}=\nu_{3d}=\nu_{3e}=\nu_{6x}=0 \quad \text{when}\  k=6,\nonumber\end{align}
where $6x$ denotes  one of the following conjugacy classes: $6a$, $6b$,  $6c$.

Let  $u\in V(\mathbb ZG)$ be a non-trivial involution. According to
(2) we get   $ \nu_{2a}+\nu_{2b}+\nu_{2c}=1$ and  by (1)
\[
\begin{split}
\mu_{0}(u,\chi_{2}) =& \frac{1}{2} (\nu_{2a}
-\nu_{2b} -\nu_{2c} + 1 );\quad \mu_{1}(u,\chi_{11}) = \frac{1}{2} (\nu_{2a} -3\nu_{2b} +\nu_{2c} + 3 );\\
\mu_{1}(u,\chi_{2}) = & \frac{1}{2} (1 -\nu_{2a}
+\nu_{2b} +\nu_{2c} );\quad\mu_{0}(u,\chi_{11}) = \frac{1}{2} (3 -\nu_{2a} +3\nu_{2b} -\nu_{2c});\\
\mu_{0}(u,\chi_{7}) = &\nu_{2a} + 1;\qquad\qquad \mu_{1}(u,\chi_{10}) =  \frac{1}{2} (\nu_{2a} + 3\nu_{2b} -\nu_{2c} + 3 ); \\
\mu_{1}(u,\chi_{7}) =&1-\nu_{2a};\qquad\qquad
\mu_{0}(u,\chi_{10}) =  \frac{1}{2} ( -\nu_{2a}-3\nu_{2b} +\nu_{2c} + 3 ).\\
\end{split}
\]
Since $\mu_{i}(u,\chi_{j})\geq 0$,  it follows that
$$
(\nu_{2a},\nu_{2b},\nu_{2c})\in \{(0, 0, 1),\ (0, 1, 0 ),\ (1, 0, 0 ) \}. $$

Let  $u\in V(\mathbb ZG)$ be a non-trivial unit of order $3$. Put
$$
\nu_1=\nu_{3a},\quad  \nu_2=\nu_{3b}, \quad \nu_3=\nu_{3c},\quad
\nu_4=\nu_{3d}, \quad \nu_5=\nu_{3e}. $$
According to (2) we get $\sum_{j=1}^5\nu_j=1$. By (1)
\[
\begin{split}
\mu_{0}(u,\chi_{3}) = & \frac{1}{3} ( - \nu_1 - \nu_2 + 2 \nu_3 - \nu_4 - \nu_5 + 1 ) \geq 0;\\
\mu_{1}(u,\chi_{3}) = &\frac{1}{3} ( - \nu_1 + 2 \nu_2 - \nu_3 - \nu_4 + 2 \nu_5 + 1 ) \geq 0;\\
\mu_{2}(u,\chi_{3}) = &\frac{1}{3} ( 2 \nu_1 - \nu_2 - \nu_3 + 2 \nu_4 - \nu_5 + 1 ) \geq 0;\\
\mu_{0}(u,\chi_{7}) = &\frac{1}{3} ( 4 \nu_1 + 4 \nu_2 -2 \nu_3 -2 \nu_4 -2 \nu_5 + 2 ) \geq 0;\\
\mu_{1}(u,\chi_{7}) = &\frac{1}{3} ( -2 \nu_1 -2 \nu_2 +  \nu_3 +  \nu_4 +  \nu_5 + 2 ) \geq 0;\\
\mu_{0}(u,\chi_{8}) = &\frac{1}{3} ( -2 \nu_1 -2 \nu_2 -2 \nu_3 +  \nu_4 +  \nu_5 + 2 ) \geq 0;\\
\mu_{1}(u,\chi_{8}) = &\frac{1}{3} ( -2 \nu_1 + 4 \nu_2 +  \nu_3 +  \nu_4 -2 \nu_5 + 2 ) \geq 0;\\
\mu_{2}(u,\chi_{8}) = &\frac{1}{3} ( 4 \nu_1 -2 \nu_2 +  \nu_3 -2 \nu_4 +  \nu_5 + 2 ) \geq 0;\\
\mu_{0}(u,\chi_{10}) = &2 \nu_3 + 1 \geq 0;\qquad \qquad \mu_{1}(u,\chi_{10}) = - \nu_3+ 1  \geq 0,\\
\end{split}
\]
so we obtain only the following trivial solution:
\[
\begin{split}
(\nu_1, \nu_2, \nu_3, \nu_4, \nu_5)\in  \{ &( 0, 0, 0, 0, 1 ),\;
( 0, 0, 0, 1, 0), \;  (0, 0, 1, 0, 0),\\
&(0, 1, 0, 0, 0),\; (1, 0, 0, 0, 0 ) \}. \end{split} \]
Let  $u\in V(\mathbb ZG)$ be a non-trivial unit of order $6$. Clearly,
$\chi(u^3)$ coincides either with
  $\chi(2a)$ or $\chi(2b)$ or $\chi(2c)$ and
$\chi(u^2)$ coincides  either with    $\chi(3a)$ or $\chi(3b)$ or $\chi(3c)$ or
$\chi(3d)$ or $\chi(3e)$. By (1) we obtain $15$ systems of inequalities. These have
no integral solutions, except for the case $\chi(u^3)=\chi(2b)$ and
$\chi(u^2)=\chi(3a)$. In this exceptional case, it is easy to see that
\[
\begin{split}
\mu_1(u,\chi_2)&= \mu_2(u,\chi_2),\quad  \mu_0(u,\chi_3)=
\mu_3(u,\chi_3),\quad \mu_1(u,\chi_3)= \mu_4(u,\chi_3),\\
\mu_0(u,\chi_4)&= \mu_3(u,\chi_4),\quad  \mu_1(u,\chi_4)=
\mu_4(u,\chi_4),\quad  \mu_1(u,\chi_7)= \mu_2(u,\chi_7),\\
\mu_0(u,\chi_8)&= \mu_3(u,\chi_8), \quad \mu_1(u,\chi_8)=
\mu_4(u,\chi_8).
\end{split}
\]
Using these additional relations, we
obtain that
$$
(\nu_{2a},\nu_{3a}, \nu_{2c}, \nu_{2d},\nu_{4a},\nu_{6a},\nu_{2e},
\nu_{4b})=( 0, 0, 0, 0, 0, 1, 0, 0, 0, 0, 0).
$$
Thus Theorem D completes the proof.
\end{proof}

\section{The octahedral group $S_4 \times  C_2$}

In this section we want to give an example of a group for which the
Luthar--Passi method is not sufficient to proof ZC--1.

The group $G=S_4 \times  C_2$ has the following character table
(easily verified with \cite{8}):

\bigskip
\centerline{ \vbox {{ \offinterlineskip \halign { $#$ &\tvon \vrule
\quad \hfil  $#$ &  \vrule   \quad \hfil  $#$ & \vrule \quad \hfil
$#$ &  \vrule   \quad \hfil  $#$  & \vrule \quad \hfil $#$ &
\vrule   \quad \hfil  $#$ & \vrule \quad \hfil $#$ & \vrule \quad
\hfil  $#$ & \vrule \quad \hfil $#$ & \vrule \quad \hfil $#$ \cr
\noalign {\hrule}  & 1a & 2a & 2b & 3a & 2c & 2d & 4a
& 6a & 2e& 4b \fvon\cr \noalign {\hrule}
\chi_1    &  1 & 1  & 1 & 1 & 1 & 1&  1&  1&  1&  1 \cr
\chi_2    &  1 & -1  & -1&  1&  1&  1 & -1 & -1 & -1&  1 \cr
\chi_3    & 1 & -1 & 1&  1&  1 & -1  & -1 & 1&  1& -1 \cr
\chi_4 & 1 & 1 & -1&  1&  1 & -1  & 1  & -1& -1&  -1 \cr
\chi_5 & 2 & 0 & 2 & -1 & 2 & 0&  0& -1&  2&  0 \cr
\chi_6 & 2 & 0  & -2 & -1 & 2 & 0 & 0&  1&  -2& 0\cr
\chi_7    & 3 & -1  & -3 & 0 & -1&  1&  1  & 0&  1&  -1 \cr
\chi_8    & 3 & -1  & 3 & 0 & -1&  -1&  1  & 0&  -1&  1 \cr
\chi_9    & 3 & 1  & -3 & 0 & -1&  -1&  -1  & 0&  1&  1 \cr
\chi_{10}    & 3 & 1  & 3 & 0 & -1&  1&  -1  & 0&  -1& -1 \cr
\noalign {\hrule} }}} }
\bigskip

According to Theorem B the values 2, 3, 4, 6 and $12$ are possible
orders of a non-trivial normalized torsion unit $u$ of $\mathbb Z G$. By Theorem D units of
order 3 are conjugate to group elements. The Luthar--Passi method gives a positive answer of ZC--1 for units of order
4, 6 and 12, provided that the conjecture holds for elements of
order $2$. For an involution $u$ of $\mathbb Z G$ the Luthar--Passi
method gives $22$ possible sets of partial augmentations of $u$,
that would be in contradiction to ZC--1. In $20$ of these cases the
sum of partial augmentations of the classes of elements of order
four is different from $0$. Looking at the reduction $G \to G/C_2$
we see that the induced surjective ring homomorphism $\mathbb Z G \to \mathbb Z[G/C_2]= \mathbb Z S_4$ maps the conjugacy classes of order $4$ of $S_4 \times C_2$ onto the unique conjugacy class $C$ of elements
of order~4 of $S_4 $. Thus the image of $u$ is a unit which has at
the conjugacy class $C$ a partial augmentation different from zero.
But for $S_4$ the conjecture  ZC--1 holds \cite{6}.

In the remaining two cases the involution $u$ has one of the
following sets of partial augmentations different from $0$:
\[
\begin{split}
 (\nu_{2a},\nu_{4a},\nu_{4b}) &= (1,-1,1); \\
 (\nu_{2d},\nu_{4a},\nu_{4b}) &= (1,1,-1). \\
\end{split}
\]
If one considers now only the reduction maps $G \to G/N $, where $N$ is
an arbitrary normal subgroup of $G$, then no contradiction occurs.
In fact an involution of this type would map modulo $C_2$ to a unit
which has the same type as a transposition of $\mathbb Z S_4 $. With
respect to the other reductions the image has also always the same
type as a trivial unit.

Using the fact that ZC--1 is valid for $S_4$ \cite{6}, the
conjecture ZC--1 follows however for $S_4 \times C_2$ from the following
result \cite{9}.

\begin{Proposition}\label{Prop2}
The conjecture ZC--1 holds for
$G \times C_2$ provided it holds for $G $.
\end{Proposition}

\begin{proof}
Let $C_2 = \gp{t} $. Let $\kappa $ be the surjective
ring homomorphism $\mathbb Z[G \times C_2]\longrightarrow \mathbb Z G $
induced by the projection $\pi: G \times C_2 \longrightarrow G $
and denote by $\alpha $     the ring homomorphism which is induced by the identity
on $G$ and $t \mapsto -1 \in \mathbb Z G $.

Two conjugacy classes of $G \times C_2$ have the same image
$\overline{C}$ in $G$ if and only if they are of the form $C_1 = C \times \{1\}$ and $C_2 = C \times \{t\} $, where $C$ is a conjugacy class of $G $. Let $u \in V(\mathbb Z G)$ be of order $k $. Denote the
partial augmentations of $u$ with respect to $C_1$ and $C_2$ by
$\nu_1$ and $\nu_2 $. The image of $u$ under $\kappa $ has with
respect to $\overline{C}$ the partial augmentation $\nu_1 + \nu_2 $
and under $\alpha $ the partial augmentation $\nu_1 - \nu_2 $.

Because ZC--1 is valid for $\mathbb Z G$ we get that $\nu_1 + \nu_2 \in \{0,1\} $ and that $\nu_1 -\nu_2\in \{0,1,-1\}$.
Note that the unit $\alpha (u) $ may be not normalized. It follows that $\nu_1+\nu_2=0$ if and only if $\nu_1=\nu_2=0 $
and  $\nu_1 + \nu_2 \neq 0 $ if and only if $\nu_1 = 1$ and $\nu_2 = 0 $ or $\nu_1 = 0$ and $ \nu_2 = 1 $. Thus there is only one partial augmentation of $u$ which is different from zero.
By Theorem C we conclude that $u$ is conjugate to an element of $G \times C_2 $.
\end{proof}

We remark that it seems to be unknown whether the conjecture ZC--1 is valid for
a direct product $ G \times H $ if it holds for $G$ and $H $.
\newpage

\section{3-dimensional crystallographic point groups}

\begin{theore}
ZC--1 is valid for all finite 3-dimensional crystallographic point groups.
\end{theore}

\begin{proof}
By the classification of 3-dimensional
crystallographic point groups the maximal ones are $D_6 \times
C_2$, $D_4 \times C_2$, $S_4 \times C_2 $, where $D_n$ denotes the
dihedral group of order $2n$.

ZC--1 is known for $D_6 \times C_2$ and all its subgroups by \cite{19},  for the
\linebreak $2$-group $ D_4 \times C_2 $ and all its subgroups by \cite{21}. The
subgroups of $S_4 \times C_2$ are of type $G$ resp. $G \times C_2 ,$ where
$G$ is either metacyclic, or isomorphic to $A_4 $ or $S_4 .$ ZC--1 for $A_4$ has been proved in
\cite{1}, and for $S_4$ in \cite{6}. So it remains the case of the octahedral group $S_4 \times C_2 .$ Thus the
previous section completes the proof.
\end{proof}

\begin{rema}
Let $H$ be a subgroup of the finite group $G $. Suppose that ZC--1 holds for $G $. Then we may consider $\mathbb Z
H$ as subring of $\mathbb Z G$ via the inclusion of $H$ in $G $. Because ZC--1 holds, a torsion unit $u$ of
$V(\mathbb Z H)$ is conjugate within $\mathbb Q G$ to $g \in G $. If $C$ is a conjugacy
class of $G$ with $C \cap H =\varnothing$ then the partial augmentation $\varepsilon_C (u) = 0 $. Therefore $u$ is conjugate to $h \in H$ within $\mathbb Q G$. This leads naturally to the question whether this conjugation may be realized in $\mathbb Q H $. If this is the case ZC--1 would follow for $H $ and the proof of ZC--1 for all finite groups would be reduced to that one for finite symmetric groups.

In the context of the proof of the theorem above it would then suffice to establish ZC--1 for the maximal finite crystallographic point groups.

Thus it might be worthwhile to consider the following version of the Zassenhaus conjecture ZC--1:

\begin{itemize}
\item[]
\begin{itemize}
\item[{\bf ZC--$1_R$}] Suppose that $u$ is a torsion unit of $V(\mathbb Z G)$ then there exists a ring $R$ containing $\mathbb Z G$ as subring such that $u$ is conjugate in $R$ to an element of $G $.
\end{itemize}
\end{itemize}
\end{rema}

\section{Generic Character Tables}

Let $p$ be an odd prime $\geq 5.$ The group $PSL(2,p)$ has the following generic
character tables (\cite{7} \S 10, \cite{4}). \pagebreak In the second row we list the orders of
the elements in the respective conjugacy classes.

Define $\varepsilon \in\{\pm 1 \}$ by $p\equiv \varepsilon \pmod 4$.

\bigskip
\small{
\centerline{ \vbox {{ \offinterlineskip \halign { $#$ & \vrule
\quad \hfil  $#$ \;  &  \vrule   \quad \hfil  $#$\; & \vrule
\quad \hfil  $#$ \; &  \vrule   \quad \hfil  $#$ \; & \vrule
\quad \hfil  $#$\tvon \cr \noalign
{\hrule} & C_1 & C_2 & C_3 & C_4(j)  & C_5(j)\fvon  \cr \noalign {\hrule}
& 1 & p & p & \frac{p-1}{2} & \frac{p+1}{2}  \fvon \cr \noalign {\hrule}
\chi_1    &  1 & 1  & 1 & 1 & 1 \cr
\chi_2    &  p & 0  & 0  &  1&  -1 \cr
\chi_3    & \frac{p+\varepsilon}{2} & \frac{1}{2}(\varepsilon-\sqrt{\varepsilon p})   &
\frac{1}{2}(\varepsilon+\sqrt{\varepsilon p}) &  (-1)^j\frac{1+\varepsilon}{2} &  (-1)^j\frac{1-\varepsilon}{2} \cr
\chi_4 &   \frac{p+\varepsilon}{2}  & \frac{1}{2}(\varepsilon+\sqrt{\varepsilon p})   &
\frac{1}{2}(\varepsilon-\sqrt{\varepsilon p}) &  (-1)^j\frac{1+\varepsilon}{2} &  (-1)^j\frac{1-\varepsilon}{2} \cr
\chi_5(k) & p+1 & 1 & 1 & a^{2jk}+a^{-2jk}    & 0 \cr
\chi_6(k) & p-1 & -1  & -1 & 0 & -b^{2jk}-b^{-2jk} \cr
\noalign {\hrule} }}} }}%
\normalsize
\bigskip

\noindent with $ a = e^{\frac{2\pi i}{p-1}}$ and  $b = e^{\frac{2\pi i}{p+1}}$.

For the sequel let $z =  e^{\frac{2\pi i}{p}} $. By Theorem C we
know that the partial augmentations $\nu_i $ of a unit $u \in
\mathbb Z[PSL(2,p)]$ of order $p$ are zero for the classes $C_1$,
$C_4(j)$ and $C_5(j)$. Thus for an irreducible character $\chi $,
using $\nu_2 + \nu_3 = 1 $, we get $\chi (u) = \nu_2 \chi (C_2) +
(1 - \nu_2) \chi (C_3)$. Hence
$$ \chi_3 (u) = \frac{\varepsilon}{2} + (-2\nu_2 + 1) \frac{\sqrt{\varepsilon p}}{2} .$$

For the calculation of the multiplicities we use the following.

\noindent For any integer $l$,
\[
\begin{array}{cl}
Tr_{\mathbb Q (z)/\mathbb Q} \big(\sqrt{\varepsilon p} \cdot z^{-l}\big)  = \sum
\limits_{\nu=1}^{p-1 }\sqrt{\varepsilon p}\left( \displaystyle \frac{\nu}{p} \right) z^{-l \nu}\\
  = \sqrt{\varepsilon p} \, \sum\limits_{\nu=1}^{p-1 }\left(\displaystyle
 \frac{\nu}{p}\right)(z^{\nu})^{-l}
\end{array}
\]
where $\left(\frac{\nu}{p}\right)$ denotes the Legendre symbol. Thus
we get
\[
Tr_{\mathbb Q (z)/\mathbb Q} (\sqrt{\varepsilon p}\,)=0 \quad \text{for } l\equiv 0 \pmod p.
\]
{\sc C. F. Gauss}\footnote{We thank the referee for this reference.} proved in \cite{9a} that
\[
\sum\limits_{\nu=1}^{p-1}\left(\frac{\nu}{p}\right)e^{\frac{2\pi i \nu} p}=
\sum\limits_{\nu=1}^{p-1}\left(\frac{\nu}{p}\right)z^{\nu}=\sqrt{\varepsilon p}.
\]
Hence
\[\bdshs
Tr_{\mathbb Q (z)/\mathbb Q} \big(\sqrt{\varepsilon p} \cdot z^{-l}\big)=\left (\frac{-l}{p}
  \right)
\varepsilon p= \left (\frac{-l}{p} \right )p \quad \text{for } l \not\equiv 0 \pmod p.
\]
For the multiplicities we get by Theorem A
\[
\begin{tabular}{ccll}
(i) & $\mu_0(u,\chi_3)$& $= \frac{1+\varepsilon}{2}$.\\[2mm]
(ii) & $\mu_l(u,\chi_3)$& $= 1 - \nu_2$& \; if \qquad $Tr_{\mathbb Q
(z)/\mathbb Q} \big(\sqrt{\varepsilon p} \cdot z^{-l}\big) = +p.$ \\[2mm]
(iii) & $\mu_l(u,\chi_3)$& $= \nu_2$& \; if \qquad $Tr_{\mathbb Q
(z)/\mathbb Q} \big(\sqrt{\varepsilon p} \cdot z^{-l}\big) = -p .$\\[2mm]
\end{tabular}
\]
Because the multiplicities have to be non-negative it follows that $\nu_2 = 0$ or
that $\nu_2 = 1 .$ This shows that elements of order $p$ are conjugate to an element
of a group basis in $\mathbb Z [PSL(2,p)]$.

The order of a finite subgroup of $V(\mathbb Z G)$ has to divide the order of $G .$
Thus finite $p$ - subgroups of the normalized unit subgroup of $\mathbb Z[PSL(2,p)]$
have order $p .$ Summarizing we get

\begin{Proposition}\label{Prop 3}
The variation p-ZC--3 is valid for $PSL(2,p)$.
\end{Proposition}

Finally we consider p-ZC--3 in the spirit of the variation ZC-1$_R$  for permutation
groups of prime degree.

\begin{Proposition}\label{Prop 4}
Let $p$ be a prime and let $G$ be a primitive permutation group of degree $p $. Consider $G$ as subgroup of the
symmetric group $S_p$. Let $u$ be a torsion unit of order $p$ of
$V(\mathbb Z G) $. Then $u$ is conjugate within $\mathbb Q S_p$ to an
element of $G $. In particular p-ZC--3  is valid for $S_p $.
\end{Proposition}

\begin{proof}
By the remark after the Theorem it suffices to
establish the result for $G = S_p .$  The symmetric group $S_p$ has only one
conjugacy class $K$ of elements of order $p .$ By Theorem C all partial
augmentations $\varepsilon_C = 0$ for all classes $C \neq K .$ By Theorem D we get
that $u$ is conjugate to an element of $S_p $ within $\mathbb Q G .$ Because $p^2$ does
not divide $p!$, it follows that \linebreak p-ZC--3 is valid for $S_p .$
\end{proof}

\end{document}